\documentclass[preprint,12pt]{elsarticle}

\usepackage{amssymb}
\usepackage{amsmath}
\usepackage{amsthm}
\usepackage[colorlinks]{hyperref}
\AtBeginDocument{%
  \hypersetup{%
    linkcolor=blue,%
    citecolor=green,%
  }%
}
\usepackage{cleveref}
\usepackage{geometry}
\usepackage{amsthm}
\usepackage{upgreek}
\usepackage[pagewise,mathlines]{lineno}
\newcommand*\patchAmsMathEnvironmentForLineno[1]{%
  \expandafter\let\csname old#1\expandafter\endcsname\csname #1\endcsname
  \expandafter\let\csname oldend#1\expandafter\endcsname\csname end#1\endcsname
  \renewenvironment{#1}%
     {\linenomath\csname old#1\endcsname}%
     {\csname oldend#1\endcsname\endlinenomath}}%
\newcommand*\patchBothAmsMathEnvironmentsForLineno[1]{%
  \patchAmsMathEnvironmentForLineno{#1}%
  \patchAmsMathEnvironmentForLineno{#1*}}%
\AtBeginDocument{%
\patchBothAmsMathEnvironmentsForLineno{equation}%
\patchBothAmsMathEnvironmentsForLineno{align}%
\patchBothAmsMathEnvironmentsForLineno{flalign}%
\patchBothAmsMathEnvironmentsForLineno{alignat}%
\patchBothAmsMathEnvironmentsForLineno{gather}%
\patchBothAmsMathEnvironmentsForLineno{multline}%
}

\crefname{equation}{}{}

\newtheorem{theorem}{Theorem}[section]
\newtheorem{lemma}[theorem]{Lemma}

\theoremstyle{definition}
\newtheorem{definition}[theorem]{Definition}

\theoremstyle{remark}
\newtheorem{remark}[theorem]{Remark}

\numberwithin{equation}{section}

\journal{~}

\begin{document}

\begin{frontmatter}



\title{Boundary Lipschitz Regularity and the Hopf Lemma for Fully Nonlinear Elliptic Equations \tnoteref{t1}}

\author[rvt]{Yuanyuan Lian}
\ead{lianyuanyuan@sjtu.edu.cn;lianyuanyuan.hthk@gmail.com}
\author[rvt]{Kai Zhang\corref{cor1}}
\ead{zhangkaizfz@gmail.com}
\tnotetext[t1]{This research is supported by the China Postdoctoral Science Foundation (Grant No.
2021M692086 and 2022M712081), the National Natural Science Foundation of China (Grant No. 12031012, 11831003 and 12171299) and  the Institute of Modern Analysis-A Frontier Research Center of Shanghai.}

\cortext[cor1]{Corresponding author. ORCID: \href{https://orcid.org/0000-0002-1896-3206}{0000-0002-1896-3206}}

\address[rvt]{School of Mathematical Sciences, Shanghai Jiao Tong University, Shanghai, China}

\begin{abstract}
In this paper, we study the boundary regularity for viscosity solutions of fully nonlinear elliptic equations. We use a unified, simple method to prove that if the domain $\Omega$ satisfies the exterior $C^{1,\mathrm{Dini}}$ condition at $x_0\in \partial \Omega$, the solution is Lipschitz continuous at $x_0$; if $\Omega$ satisfies the interior $C^{1,\mathrm{Dini}}$ condition at $x_0$, the Hopf lemma holds at $x_0$. The key idea is that the curved boundaries are regarded as perturbations of a hyperplane. Moreover, we show that the $C^{1,\mathrm{Dini}}$ conditions are optimal.
\end{abstract}


\end{frontmatter}

\section{Introduction}\label{S1}
It is well known that the geometric properties of domains have significant influence on the boundary regularity of solutions. For elliptic equations in nondivergence form, if the domain $\Omega$ satisfies the exterior sphere condition at $x_0\in \partial \Omega$, the solution is Lipschitz continuous at $x_0$; if $\Omega$ satisfies the interior sphere condition at $x_0$, the Hopf lemma holds at $x_0$. These can be proved easily by constructing proper barriers (see \cite[P. 27 and Lemma 3.4]{MR1814364} and \cite[Lemma 1.1 and Lemma 1.2]{Safonov2008}).

The sphere condition is not the optimal geometrical condition and has been generalized. In particular, under the exterior and the interior $C^{1,\mathrm{Dini}}$ conditions (see \Cref{d-re} and \Cref{d-H-re}), Safonov \cite{Safonov2008} proved the boundary Lipschitz regularity and the Hopf lemma respectively. However, he proved for classical solutions of linear elliptic equations and the boundary Harnack inequality was applied as a main tool. Huang, Li and Wang \cite{MR3167627} also obtained the boundary Lipschitz regularity for linear elliptic equations under the exterior $C^{1,\mathrm{Dini}}$ condition. They used an auxiliary function and the iteration technique, without the usage of the boundary Harnack inequality. Lieberman \cite{MR779924} proved the Hopf lemma for linear elliptic equations under the interior $C^{1,\mathrm{Dini}}$ condition by applying the regularized distance. Narazov \cite{MR2888295} extended to equations involving the gradient in a half ball and he remarked that the results still hold for $C^{1,\mathrm{Dini}}$ conditions with the aid of the technique in \cite{MR779924}. Recently, Apushkinskaya and Nazarov \cite{MR4030346} derived the Hopf lemma for equations with lower order terms in divergence form under the interior $C^{1,\mathrm{Dini}}$ condition. For more details on the Hopf lemma and related results, we refer to the comprehensive survey \cite{MR4461367}.

In regard to the optimality of the $C^{1,\mathrm{Dini}}$ condition, Safonov \cite{Safonov2008} stated the optimality without a proof. Recently, Safonov \cite{MR3802819} proved the optimality of the interior $C^{1,\mathrm{Dini}}$ condition for the Hopf lemma. In addition, Apushkinskaya and Nazarov \cite{MR3513140} also showed the optimality of the interior $C^{1,\mathrm{Dini}}$ condition under an additional convexity assumption on the domain.

It is interesting to note that the boundary Lipschitz regularity needs a geometrical condition from the exterior but the Hopf lemma needs a geometrical condition from the interior. Moreover, both results need the same geometrical condition (sphere condition or $C^{1,\mathrm{Dini}}$ condition).

In this paper, we consider fully nonlinear elliptic equations and prove the boundary Lipschitz regularity and the Hopf lemma for viscosity solutions under the $C^{1,\mathrm{Dini}}$ conditions. These two results are proved by a unified method and the proofs are simple. Furthermore, we show that the $C^{1,\mathrm{Dini}}$ conditions are optimal with short proofs.

The key idea for the boundary Lipschitz regularity and the Hopf lemma is that the curved boundaries are regarded perturbations of a flat boundary (i.e., a hyperplane). Based on the boundary $C^{1,\alpha}$ regularity for solutions with flat boundaries (see \cite[Theorem 9.31]{MR1814364}, \cite[Theorem 4.28]{MR787227} and \cite[Lemma 7.1]{MR2254613}), the boundary regularity for solutions with curved boundaries can be obtained through a perturbation argument. The derivatives maybe enlarged (shrunk) when the boundary is curved toward the exterior (interior). If the boundary is not curved toward the exterior too much, the derivatives remain bounded. This is just the boundary Lipschitz regularity. Similarly, if the boundary is not curved toward the interior too much, the positive derivative doesn't vanish. This is just the Hopf lemma. The exterior and interior $C^{1,\mathrm{Dini}}$ conditions can guarantee that the boundary is not curved too much and they are the optimal geometrical conditions. It explains why the boundary Lipschitz regularity and the Hopf lemma require the same $C^{1,\mathrm{Dini}}$ condition, and one requires from the exterior and the other requires from the interior.

We use solutions with flat boundaries to approximate the solution and the error between them can be estimated by maximum principles. This basic perturbation idea is inspired originally by \cite{MR1351007}. The application to boundary regularity is inspired by \cite{MR3780142}.

Before the statement of our main results, we introduce some notations and notions. Denote $B_r(x_0)$ the open ball in $\mathbb{R}^n$ with center $x_0$ and radius $r$. Let $B_r=B_r(0)$, $B_r^+=B_r\cap \left\{x=(x_1,\cdots,x_n)\in \mathbb{R}^n|x_n>0\right\}$, $T_r=B_r\cap \left\{x\in \mathbb{R}^n|x_n=0\right\}$ and $T_r(x_0)=T_r+x_0$.

\begin{definition}\label{d-Dini}
A function $\omega:[0,+\infty)\rightarrow [0,+\infty)$ is called a Dini function if $\omega$ is nondecreasing and satisfies the following Dini condition for some $r_0>0$,
\begin{equation}\label{e-dini}
    \int_{0}^{r_0}\frac{\omega(r)}{r}dr<\infty.
\end{equation}
\end{definition}

Now, we give the definitions of the geometrical conditions, under which we prove our main results.
\begin{definition}[\textbf{exterior $C^{1,\mathrm{Dini}}$ condition}]\label{d-re}
We say that $\Omega$ satisfies the exterior $C^{1,\mathrm{Dini}}$ condition at $x_0\in \partial \Omega$ if there exists $r_0>0$ and a coordinate system $\{x_1,...,x_n \}$ such that $x_0=0$ in this coordinate system and
\begin{equation}\label{e-re}
B_{r_0} \cap \{x_n \leq -|x'|\omega(|x'|)\} \subset B_{r_0}\cap \Omega^c,
\end{equation}
where $\omega$ is a Dini function.
\end{definition}

\begin{definition}[\textbf{anti-exterior $C^{1,\mathrm{Dini}}$ condition}]\label{d-re-2}
We say that $\Omega$ satisfies the anti-exterior $C^{1,\mathrm{Dini}}$ condition at $x_0\in \partial \Omega$ if there exists $r_0>0$ and a coordinate system $\{x_1,...,x_n \}$ such that $x_0=0$ in this coordinate system and
\begin{equation}\label{e-re-2}
B_{r_0} \cap \{x_n >-|x'|\omega(|x'|)\} \subset B_{r_0}\cap \Omega,
\end{equation}
where $\omega(r)\rightarrow 0$ as $r\rightarrow 0$ and doesn't satisfy \cref{e-dini}, i.e., the integral is infinity.
\end{definition}

\begin{definition}[\textbf{interior $C^{1,\mathrm{Dini}}$ condition}]\label{d-H-re}
We say that $\Omega$ satisfies the interior $C^{1,\mathrm{Dini}}$ condition at $x_0\in \partial \Omega$ if there exists $r_0>0$ and a coordinate system $\{x_1,...,x_n \}$ such that $x_0=0$ in this coordinate system and
\begin{equation}\label{e-H-re}
B_{r_0} \cap \{x_n >|x'|\omega(|x'|)\} \subset B_{r_0}\cap \Omega,
\end{equation}
where $\omega$ is a Dini function.
\end{definition}

\begin{definition}[\textbf{anti-interior $C^{1,\mathrm{Dini}}$ condition}]\label{d-H-re-2}
We say that $\Omega$ satisfies the anti-interior $C^{1,\mathrm{Dini}}$ condition at $x_0\in \partial \Omega$ if there exists $r_0>0$ and a coordinate system $\{x_1,...,x_n \}$ such that $x_0=0$ in this coordinate system and
\begin{equation}\label{e-H-re-2}
B_{r_0} \cap \{x_n \leq |x'|\omega(|x'|)\} \subset B_{r_0}\cap \Omega^c,
\end{equation}
where $\omega(r)\rightarrow 0$ as $r\rightarrow 0$ and doesn't satisfy \cref{e-dini}.
\end{definition}

\begin{remark}\label{r-5}
In this paper, if we say that $\Omega$ satisfies any of \Crefrange{d-re}{d-H-re-2}, we always denote the corresponding function by $\omega_{\Omega}$ and assume that $r_0=1,\omega_{\Omega}(1)<1/4$.
\end{remark}

Next, we introduce some definitions to describe the pointwise regularity for solutions.
\begin{definition}\label{d-f}
Let $\Omega \subset \mathbb{R}^{n}$ be a bounded domain and $f$ be a function defined on $\bar{\Omega}$. We say that $f$ is Lipschitz continuous at $x_0\in \bar{\Omega}$ ($f$ is $C^{0,1}$ at $x_0$ or $f\in C^{0,1}(x_0)$) if there exists a constant $C$ such that
\begin{equation*}\label{holder}
  |f(x)-f(x_0)|\leq C|x-x_0|,~~\forall~x\in \bar{\Omega}.
\end{equation*}
Then, define $[f]_{C^{0,1}(x_0)}=C$ and $\|f\|_{C^{0,1}(x_0)}=\|f\|_{L^{\infty}(\Omega)}+[f]_{C^{0,1}(x_0)}$.

Similarly, we call that $f$ is $C^{1,\mathrm{Dini}}$ at $x_0$ or $f\in C^{1,\mathrm{Dini}}(x_0)$ if there exist a vector $l$ and a constant $C$ such that
\begin{equation*}
  |f(x)-f(x_0)-l\cdot (x-x_0)|\leq C|x-x_0|\omega(|x-x_0|),~~\forall~x\in \bar{\Omega},
\end{equation*}
where $\omega$ is a Dini function and its Dini integral is equal to $1$. Then we denote $l$ by $\nabla f(x_0)$. We define $[f]_{C^{1,\mathrm{Dini}}(x_0)}=C$ and $\|f\|_{C^{1,\mathrm{Dini}}(x_0)}=\|f\|_{L^{\infty}(\Omega)}+|l|+[f]_{C^{1,\mathrm{Dini}}(x_0)}$.

Finally, we define that $f$ is $C^{-1,\mathrm{Dini}}$ at $x_0$ or $f\in C^{-1,\mathrm{Dini}}(x_0)$ if there exists a constant $C$ such that
\begin{equation*}
\|f\|_{L^{n}(B(x_0,r))}\leq C\omega(r),~~\forall~0<r<\mathrm{diam}(\Omega),
\end{equation*}
where $\mathrm{diam}(\Omega)$ denotes the diameter of $\Omega$ and $\omega$ is a Dini function with its Dini integral equal to $1$. Then we define $\|f\|_{C^{-1,\mathrm{Dini}}(x_0)}=C$.
\end{definition}

\begin{remark}\label{r-6}
In this paper, if we say that a function $f$ is $C^{1,\mathrm{Dini}}$ or $C^{-1,\mathrm{Dini}}$ at $x_0$, we use $\omega_f$ to denote the corresponding Dini function.
\end{remark}

Since we consider viscosity solutions of fully nonlinear elliptic equations, we briefly introduces some notions in this respect (see \cite{MR1376656,MR1351007,MR1118699} for more details). We call $F:S^n\rightarrow \mathbb{R}$ a fully nonlinear uniformly elliptic operator with ellipticity constants $0<\lambda\leq \Lambda$ if
\begin{equation*}
    \lambda\|N\|\leq F(M+N)-F(M)\leq \Lambda\|N\|,~~~~\forall~M,N\in S^n,~N\geq 0,
\end{equation*}
where $S^n$ denotes the set of $n\times n$ symmetric matrices; $\|N\|$ is the spectral radius of $N$ and $N\geq 0$ means the nonnegativeness. A constant is called universal if it depends only on the dimension $n$ and the ellipticity constants $\lambda$ and $\Lambda$.

Next, we introduce the notion of viscosity solution.
\begin{definition}\label{d-viscoF}
Let $u\in C(\Omega)$ and $f\in L^{n}(\Omega)$. We say that $u$ is an $L^n$-viscosity subsolution (resp., supersolution) of
\begin{equation}\label{FNE}
F(D^2u)=f ~~~~\mbox{in}~\Omega,
\end{equation}
if
\begin{equation*}
  \begin{aligned}
    \mathrm{ess}~\underset{y\to x}{\lim\sup}\left(F(D^2\varphi(y))-f(y)\right)\geq 0\\
    \left(\mathrm{resp.},~\mathrm{ess}~\underset{y\to x}{\lim\inf}\left(F(D^2\varphi(y))-f(y)\right)\leq 0\right)
  \end{aligned}
\end{equation*}
provided that for $\varphi\in W^{2,n}(\Omega)$, $u-\varphi$ attains its local maximum (resp., minimum) at $x\in\Omega$.

We call $u$ an $L^n$-viscosity solution of \cref{FNE} if it is both an $L^n$-viscosity subsolution and supersolution of \cref{FNE}.
\end{definition}

Next, we give the definitions of the Pucci's extremal operators and the Pucci's class.
\begin{definition}\label{d-Sf}
For $M\in S^n$, denote its eigenvalues by $\lambda_i$ ($1\leq i\leq n$) and define the Pucci's extremal operators:
\begin{equation*}
\begin{aligned}
  &\mathcal{M}^+(M,\lambda,\Lambda)=\Lambda\left(\sum_{\lambda_i>0}\lambda_i\right)
+\lambda\left(\sum_{\lambda_i<0}\lambda_i\right), \\
  &\mathcal{M}^-(M,\lambda,\Lambda)=\lambda\left(\sum_{\lambda_i>0}\lambda_i\right)
+\Lambda\left(\sum_{\lambda_i<0}\lambda_i\right).
\end{aligned}
\end{equation*}

The Pucci's class is defined as follows. We say that $u\in \underline{S}(\lambda,\Lambda,f)$ if $u$ is an $L^n$-viscosity subsolution of
\begin{equation*}
  \mathcal{M}^+(D^2u,\lambda,\Lambda)= f.
\end{equation*}
Similarly, we denote $u\in \bar{S}(\lambda,\Lambda,f)$ if $u$ is an $L^n$-viscosity supersolution of
\begin{equation*}
  \mathcal{M}^-(D^2u,\lambda,\Lambda)= f.
\end{equation*}
We also define
\begin{equation}\label{SC2Sf}
\begin{aligned}
  &S(\lambda,\Lambda,f)=\underline{S}(\lambda,\Lambda,f)\cap \bar{S}(\lambda,\Lambda,f),\\
  &S^*(\lambda,\Lambda,f)=\underline{S}(\lambda,\Lambda,-|f|)\cap \bar{S}(\lambda,\Lambda,|f|).
\end{aligned}
\end{equation}

Clearly, $S(\lambda,\Lambda,f)\subset S^*(\lambda,\Lambda,f)$. The importance of the Pucci's class lies in that any viscosity solution $u$ of \cref{FNE} belongs to $S(\lambda,\Lambda,f)$.
\end{definition}

Now, we state our main results and the first is the boundary Lipschitz regularity:
\begin{theorem}[\textbf{Boundary Lipschitz regularity}]\label{t-2}
Suppose that $\Omega$ satisfies the exterior $C^{1,\mathrm{Dini}}$ condition at $0\in \partial \Omega$ and $u\in C(\bar{\Omega}\cap B_1)$ satisfies
\begin{equation}
\left\{\begin{aligned}
&u\in S^*(\lambda,\Lambda,f)&& ~~\mbox{in}~~\Omega\cap B_1;\\
&u=g&& ~~\mbox{on}~~\partial \Omega\cap B_1,
\end{aligned}\right.
\end{equation}
where $g\in C^{1,\mathrm{Dini}}(0)$ and $f\in C^{-1,\mathrm{Dini}}(0)$.

Then $u$ is $C^{0,1}$ at $0$ and
\begin{equation*}
  |u(x)-u(0)|\leq C |x|\left(\|u\|_{L^{\infty }(\Omega)}+\|f\|_{C^{-1,\mathrm{Dini}}(0)}+\|g\|_{C^{1,\mathrm{Dini}}(0)}\right), ~~\forall ~x\in \Omega\cap B_{\delta},
\end{equation*}
where $C$ and $\delta$ depend only on $n, \lambda, \Lambda,\omega_f, \omega_g$ and $\omega_{\Omega}$.
\end{theorem}
\begin{remark}\label{r-1}
Note that if $f\in L^p(\Omega\cap B_1)$ with $p>n$, then $f\in C^{-1,\mathrm{Dini}}(0)$.
\end{remark}
\begin{remark}\label{r-4}
In this theorem, we don't require that $u$ satisfies an equation and it is enough that $u$ belongs to the Pucci class.
\end{remark}
\begin{remark}\label{r-9}
In this paper, we only consider the simplest fully nonlinear equations to show the idea clearly and we will study equations involving lower terms in near future.
\end{remark}

The exterior $C^{1,\mathrm{Dini}}$ condition is optimal. Indeed, we have the following result.
\begin{theorem}[\textbf{Anti-Lipschitz regularity}]\label{t-3}
Suppose that $\Omega$ satisfies the anti-exterior $C^{1,\mathrm{Dini}}$ condition at $0\in \partial \Omega$ and $u$ satisfies
\begin{equation*}
M^{-}(D^2u,\lambda,\Lambda)\leq 0 ~~\mbox{in}~~\Omega\cap B_1~~ (\mathrm{i.e.,}~~ u\in \bar{S}(\lambda,\Lambda,0))
\end{equation*}
with $u(0)=0$ and $u\geq 0$ in $\Omega\cap B_1$.

Then for any $l=(l_1,...,l_n)\in \mathbb{R}^n$ with $|l|=1$ and $l_n>0$,
\begin{equation}\label{e-L-main}
  u(rl)\geq l_n\left(\inf_{x\in B_{1/8}(e_n/2)}u(x)\right) r\tilde{\omega}(r), ~~\forall~ 0<r<\delta,
\end{equation}
where $\tilde{\omega}$ and $\delta$ depend only on $n, \lambda, \Lambda$ and $\omega_{\Omega}$; $\tilde{\omega}(r)\rightarrow +\infty$ as $r\rightarrow 0$ and $e_n=(0,\cdots,0,1)$.
\end{theorem}
\begin{remark}\label{r-7}
Suppose that $\partial \Omega\cap B_1=\left\{x\big |x_n=-|x'|\omega(|x'|)\right\}\cap B_1$. Then, if $\omega$ is a Dini function, the boundary Lipschitz regularity holds; if $\omega$ is not a Dini function, the boundary Lipschitz regularity fails. Hence, the exterior $C^{1,\mathrm{Dini}}$ condition is optimal for the boundary Lipschitz regularity. Similarly, the interior $C^{1,\mathrm{Dini}}$ condition is optimal for the Hopf lemma.
\end{remark}

For the Hopf lemma, we have
\begin{theorem}[\textbf{Hopf lemma}]\label{t-H-2}
Suppose that $\Omega$ satisfies the interior $C^{1,\mathrm{Dini}}$ condition at $0\in \partial \Omega$ and $u$ satisfies
\begin{equation*}
M^{-}(D^2u,\lambda,\Lambda)\leq 0 ~~\mbox{in}~~\Omega\cap B_1~~ (\mathrm{i.e.,}~~ u\in \bar{S}(\lambda,\Lambda,0))
\end{equation*}
with $u(0)=0$ and $u\geq 0$ in $\Omega\cap B_1$.

Then for any $l=(l_1,...,l_n)\in \mathbb{R}^n$ with $|l|=1$ and $l_n>0$,
\begin{equation}\label{e-H-main}
  u(rl)\geq cl_n \left(\inf_{x\in B_{1/8}(e_n/2)}u(x)\right)r, ~~\forall~ 0<r<\delta,
\end{equation}
where $c>0$ depends only on $n, \lambda, \Lambda$ and $\omega_{\Omega}$, and $\delta$ depends also on $l$.
\end{theorem}

\begin{theorem}[\textbf{Anti-Hopf lemma}]\label{t-H-3}
Suppose that $\Omega$ satisfies the anti-interior $C^{1,\mathrm{Dini}}$ condition at $0\in \partial \Omega$ and $u\in C(\bar{\Omega}\cap B_1)$ satisfies
\begin{equation*}
M^{+}(D^2u,\lambda,\Lambda)\geq 0 ~~\mbox{in}~~\Omega\cap B_1~~ (\mathrm{i.e.,}~~ u\in \underline{S}(\lambda,\Lambda,0))
\end{equation*}
with $u\leq 0$ on $\partial \Omega \cap B_1$.

Then
\begin{equation}\label{e-H-main-2}
u(x)\leq \sup_{\Omega\cap B_{1}} u^+\cdot |x|\tilde\omega(|x|), ~~\forall~ x\in \Omega\cap B_{\delta},
\end{equation}
where $\tilde{\omega}$ and $\delta$ depend only on $n, \lambda, \Lambda$ and $\omega_{\Omega}$, and $\tilde{\omega}(r)\rightarrow 0$ decreasingly as $r\rightarrow 0$.
\end{theorem}

\begin{remark}\label{r-2}
For linear equations, \Cref{t-3} and \Cref{t-H-3} have been pointed out by Safonov \cite{Safonov2008} and \Cref{t-H-3} was proved in \cite{MR3802819} under the following more general geometrical condition.

Assume that $\partial \Omega\cap B_1=\left\{x\big |x_n=|x'|\omega(x')\right\}\cap B_1$ where $\omega$ is a nonnegative function defined on $T_1$. For $k\geq 0$, define $V_k=T_{2^{-k}}\backslash T_{2^{-k-1}}$. Suppose that for any $k\geq 0$, there exists $T_k:=T_{\varepsilon_02^{-k}}(y_k)\subset V_k$ for some $y_k\in V_k$ and fixed $\varepsilon_0>0$ such that
\begin{equation*}
  \sum_{k=0}^{\infty} \inf_{T_k} \omega=+\infty.
\end{equation*}
Safonov proved that for linear equations, \Cref{t-H-3} holds under above geometrical condition.

We point out that our result also hold under above geometrical condition. It can be proved by the same method in this paper almost without any modification. We refer to the proof of \Cref{t-H-3} for details.
\end{remark}

\begin{remark}\label{r-8}
Apushkinskaya and Nazarov also proved \Cref{t-H-3} for strong solutions of linear equations under the following geometrical condition.

Assume that $\partial \Omega\cap B_1=\left\{x\big |x_n=F(x')\right\}\cap B_1$, where $F\geq 0$ is a convex function with $F(0)=0$. Define
\begin{equation*}
  \omega(r)=\sup_{x'\in T_r} \frac{F(x')}{|x'|}.
\end{equation*}
They proved that if $\omega$ is not a Dini function, then the Hopf lemma can't hold.

We point out that this geometrical condition is a special case of that in \Cref{r-2}. Indeed, since $F$ is convex, for any $0<r<1$, $\omega$ attains its maximum at the boundary of $T_r$. Hence, there exists $x_0'\in \partial T_{r/2}$ such that $\omega(r/2)=2F(x_0')/r$. Without loss of generality, we assume that $x_0'=(r/2,0,...,0)$. By the convexity of $F$ again,
\begin{equation*}
\left\{x'\big | F(x')\leq \frac{r}{2}\omega(\frac{r}{2})\right\}\cap T_r
\subset \left\{x'\big |x_1\leq \frac{r}{2}\right\} \cap T_r.
\end{equation*}
Hence,
\begin{equation*}
F(x')> \frac{r}{2}\omega(\frac{r}{2}) ~\mbox{in}~\left\{x'\big |x_1>\frac{r}{2}\right\} \cap T_r.
\end{equation*}

Therefore, for any $k\geq 0$, there exists $y_k\in T_{2^{-k}}\backslash T_{2^{-k-1}}$ such that
\begin{equation*}
\inf_{T_k} F\geq 2^{-k-1}\omega(2^{-k-1}),
\end{equation*}
where $T_k=T_{2^{-k-2}}(y_k)$. Thus
\begin{equation*}
\sum_{k=0}^{\infty} \inf_{T_k}\frac{F(x')}{|x'|}\geq \sum_{k=0}^{\infty} \omega(2^{-k-1})=+\infty.
\end{equation*}
That is, the geometrical condition in \Cref{r-2} holds.
\end{remark}

\begin{remark}\label{r-3}
Similar Dini conditions appear in many regularity results. Consider the following typical example:
\begin{equation*}
  \Delta u=f ~~\mbox{in}~~B_1.
\end{equation*}

It is well known that the continuity of $f$ at $0$ doesn't imply the existence of the second derivatives of $u$ at $0$ (see \cite[Problem 4.9]{MR1814364}). However, if the modulus of continuity of $f$ at $0$ is a Dini function, the second derivatives of $u$ at $0$ exist. Furthermore, if the modulus of continuity of $f$ in $B_1$ is a Dini function, $u$ belongs $C^2(\bar{B}_{1/2})$ (see \cite{MR2273802}).

Since we study partial differential equations, intuitively, a solution is obtained from a ``integration'' process. That is, in some sense, the solution is the ``integral'' of the coefficients, the right hand function and the boundary values etc. During the ``integration'' process, the boundedness of some integral or series is required. In general, the Dini condition \cref{e-dini} can guarantee the boundedness.

In the case of this paper, both the boundary Lipschitz regularity and the Hopf lemma concern the first derivatives. Hence, we need that the ``first derivatives'' satisfies the Dini condition. That is, the $C^{1,\mathrm{Dini}}$ condition is essentially necessary. From the view point of scaling, $\nabla u$ is equivalent to $\omega_g$ and $\|f\|_{L^n}$. Hence, the Dini conditions on them are also necessary for the boundary Lipschitz regularity.
\end{remark}

\section{Boundary Lipschitz regularity}

In the following two sections, we give the detailed proofs of the main results. For both of the boundary Lipschitz regularity and the Hopf lemma, we use solutions with flat boundaries (i.e., $v$ in the proofs) to approximate the solution $u$. Then the error between $u$ and $v$ (i.e., $w$ in the proofs) can be estimated by maximum principles. By an iteration argument, the boundary regularity for $u$ is obtained.

For the boundary Lipschitz regularity, the right hand function $f$, the boundary value $g$ and the curved boundary $\partial \Omega$ are regarded perturbations of $0$, $0$ and a hyperplane (see the definition of $v$ in the proof). This is inspired directly by \cite{MR3780142}. For the Hopf lemma, since the solution is nonnegative and the right hand is $0$, it is easier to prove.

On the other hand, in the cases of \Cref{t-3} and \Cref{t-H-3}, the boundaries are curved too much, which can't be approximated by flat boundaries. On the contrast, the curved boundaries provide more growth for the solution in case of \Cref{t-3} (see $v_2$ in the proof) and provide more decay in the case of \Cref{t-H-3} (see $v$ in the proof), which lead to the fails of boundary Lipschitz regularity and the Hopf lemma respectively.

Before proving the boundary regularity, we introduce some preliminary lemmas. The first states the boundary $C^{1,\alpha}$ regularity for solutions with flat boundaries. It was first proved by Krylov \cite{MR688919} and further simplified by Caffarelli (\cite[Theorem 9.31]{MR1814364} and the last paragraph of the Notes on Page 255; see also \cite[Theorem 4.28]{MR787227} and the following paragraph on Page 43). We will use the solutions in this lemma to approximate the solutions in \Cref{t-2} and \Cref{t-H-2}.
\begin{lemma}\label{l-1}
Let $u$ satisfy
\begin{equation*}
\left\{\begin{aligned}
&u\in S(\lambda,\Lambda,0)&& ~~\mbox{in}~~B_1^+;\\
&u=0&& ~~\mbox{on}~~T_1.
\end{aligned}\right.
\end{equation*}

Then $u$ is $C^{1,\alpha}$ at $0$ and there exists a constant $a$ such that
\begin{equation*}
  |u(x)-ax_n|\leq C |x|^{1+\alpha}\|u\|_{L^{\infty }(B_1^+)}, ~~\forall ~x\in B_{1/2}^+
\end{equation*}
and
\begin{equation*}
  |a|\leq C\|u\|_{L^{\infty }(B_1^+)},
\end{equation*}
where $0<\alpha<1$ and $C$ are universal constants.
\end{lemma}

The next lemma shows that the derivative is not vanishing on the boundary for a positive supersolution, which can be proved easily by constructing an appropriate barrier.
\begin{lemma}\label{le-H-1}
Let $u\geq 0$ satisfy
\begin{equation}
  \left\{\begin{aligned}
    M^{-}(D^2u,\lambda,\Lambda)&\leq 0~~&&\mbox{in}~~B_1^+;\\
    u&\geq 1~~&&\mbox{in}~~B_{1/8}(e_n/2);\\
    u&=0~~&&\mbox{on}~~T_1.
\end{aligned}\right.
\end{equation}
Then
\begin{equation*}
u(x)\geq cx_n ~~\mbox{in}~~B_{\delta}^+,
\end{equation*}
where $c>0$ and $\delta$ are universal constants.
\end{lemma}

Finally, we state two variations of \Cref{le-H-1}:
\begin{lemma}\label{le-H-2-2}
Let $u\geq 0$ satisfy
\begin{equation*}
  \left\{\begin{aligned}
    M^{-}(D^2u,\lambda,\Lambda)&\leq 0~~&&\mbox{in}~~B_1^+;\\
    u&\geq 1~~&&\mbox{on}~~T_{1/4}+e_n/2;\\
    u&=0~~&&\mbox{on}~~T_1.
\end{aligned}\right.
\end{equation*}
Then
\begin{equation*}
u(x)\geq cx_n ~~\mbox{in}~~B_{\delta}^+,
\end{equation*}
where $0<\delta<1/4$ and $c>0$ are universal constants.
\end{lemma}

\begin{lemma}\label{le-H-2}
Let $u$ satisfy
\begin{equation*}
  \left\{\begin{aligned}
    M^{-}(D^2u,\lambda,\Lambda)&\leq 0~~\mbox{in}~~B^+;\\
    u&=1~~\mbox{on}~~T;\\
    u&=0~~\mbox{on}~~\partial B^+\backslash T,
\end{aligned}\right.
\end{equation*}
where
\begin{equation*}
B^+=B_1^+\cap \left\{(x',x_n)\big| |x'|<\frac{1}{4}~\mbox{or}~x_n>a\right\}, T=T_{1}\backslash T_{1/4}+ae_n
\end{equation*}
and $0\leq a<1/2$. Then
\begin{equation*}
u(x)\geq cx_n ~~\mbox{in}~~B_{\delta}^+,
\end{equation*}
where $0<\delta<1/4$ and $c>0$ are universal constants.
\end{lemma}
We also omit the proof of this lemma. In fact, with the aid of \Cref{le-H-1}, it is enough to note that $u\geq c$ in $B_{1/8}(e_n/2)$ where $c$ is independent of $a$.

Now, we prove the main results. First, we prove the boundary Lipschitz regularity and give the~\\
\noindent\textbf{Proof of \Cref{t-2}.} Let $\omega(r)=\max \left\{\omega_{\Omega}(r),\omega_g(r),\omega_f(r)\right\}$. From the Dini condition \cref{e-dini}, there exists $r_1>0$ such that
\begin{equation}\label{e-dini-2}
 \omega(r_1)\leq c_0 ~~\mbox{and}~~ \int_{0}^{r_1}\frac{\omega(r)}{r}dr\leq c_0,
\end{equation}
where $c_0\leq 1/4$ is a small universal constant to be specified later. By a proper scaling, we assume that $r_1=1$. Furthermore, we assume that $u(0)=g(0)=0$ and $\nabla g(0)=0$. Otherwise, we may consider $v:=u-g(0)-\nabla g(0)\cdot x$, which satisfies the same equation.

Let $M=\|u\|_{L^{\infty }(\Omega)}+\|f\|_{C^{-1,\mathrm{Dini}}(0)}+[g]_{C^{1,\mathrm{Dini}}(0)}$ and $\Omega _{r}=\Omega \cap B_{r}$. To prove that $u$ is $C^{0,1}$ at $0$, we only need to prove the following: there exist universal constants $0<\alpha_{0}, \eta < 1$, $\bar{C}$ and $\hat{C}$, and a nonnegative sequence $\{a_k\}$ ($k\geq -1$) such that for all $k\geq 0$,

\begin{equation}\label{e1.16}
\sup_{\Omega _{\eta^{k}}}(u-a_kx_n)\leq \hat{C} M \eta ^{k}A_k
\end{equation}
and
\begin{equation}\label{e1.17}
|a_k-a_{k-1}|\leq \bar{C}\hat{C}MA_k,
\end{equation}
where
\begin{equation}\label{e-Ak}
A_0=c_0 ~~\mbox{and}~~A_k=\max(\omega(\eta^{k}),\eta^{\alpha_0} A_{k-1}) (k\geq 1).
\end{equation}

Indeed, from the Dini condition \cref{e-dini-2}, we have $\Sigma_{k}A_k<+\infty$. Then from \cref{e1.16} and \cref{e1.17}, there exists a nonnegative constant $a$ such that
\begin{equation*}
\sup_{\Omega _{r}}(u-ax_n)\leq C M r.
\end{equation*}
Hence,
\begin{equation*}
\sup_{\Omega _{r}}u\leq C M r.
\end{equation*}
In addition,
\begin{equation*}
\inf_{\Omega _{r}}u\geq -C M r
\end{equation*}
can be proved similarly. Therefore,
\begin{equation*}
  \|u\|_{L^{\infty}(\Omega_{r})}\leq CMr.
\end{equation*}
That is, $u$ is $C^{0, 1}$ at $0$.

Now, we prove \cref{e1.16} and \cref{e1.17} by induction. For $k=0$, by setting $a_{-1}=a_0=0$, they hold clearly provided
\begin{equation}\label{e.21}
\hat{C}c_0\geq 1.
\end{equation}
Suppose that they hold for $k$. We need to prove that they hold for $k+1$.

In the following proof, $C,C_1,C_2...,$ always denote universal constants. Let $r=\eta ^{k}/2$, $\tilde{B}^{+}_{r}=B^{+}_{r}-r\omega(r)e_n $, $\tilde{T}_r=T_r-r\omega(r) e_n$ and $\tilde{\Omega }_{r}=\Omega \cap \tilde{B}^{+}_{r}$ where $e_n=(0,0,...,0,1)$. Note that $\omega(\eta) \leq \omega(1)\leq c_0\leq 1/4$ and then $\Omega _{r/2}\subset \tilde{\Omega }_{r}\subset \Omega_{2r}$.

Let $v$ solve
\begin{equation*}
\left\{\begin{aligned}
 &M^{+}(D^2v,\lambda,\Lambda)=0 &&\mbox{in}~~\tilde{B}^{+}_{r}; \\
 &v=0 &&\mbox{on}~~\tilde{T}_{r};\\
 &v=\hat{C} M \eta ^{k}A_k &&\mbox{on}~~\partial \tilde{B}^{+}_{r}\backslash \tilde{T}_{r}.
\end{aligned}
\right.
\end{equation*}
Let $w=u-a_kx_n-v$. Then $w$ satisfies (note that $v\geq 0$ in $\tilde{B}^{+}_{r}$)
\begin{equation*}
    \left\{
    \begin{aligned}
      &w\in \underline{S}(\lambda,\Lambda , -|f|) &&\mbox{in}~~ \Omega \cap \tilde{B}^{+}_{r}; \\
      &w\leq g-a_kx_n &&\mbox{on}~~\partial \Omega \cap \tilde{B}^{+}_{r};\\
      &w\leq 0 &&\mbox{on}~~\partial \tilde{B}^{+}_{r}\cap \bar{\Omega}.
    \end{aligned}
    \right.
\end{equation*}

In the following arguments, we estimate $v$ and $w$ respectively. By the boundary $C^{1,\alpha}$ estimate for $v$ (see \Cref{l-1}) and the maximum principle, there exist a universal constant $0<\alpha<1$ and a constant $\bar{a}\geq 0$  such that
\begin{equation*}
\begin{aligned}
    \|v-\bar{a}(x_n+r\omega(r))\|_{L^{\infty }(\Omega _{2\eta r})}&\leq C\frac{(2\eta r)^{1+ \alpha}}{r^{1+ \alpha}}\|v\|_{L^{\infty }( \tilde{B}^{+}_{r})}\\
    &\leq C\eta ^{\alpha-\alpha_0 }\cdot \hat{C}M\eta ^{(k+1)}\eta^{\alpha_0}A_k\\
     &\leq C_1\eta ^{\alpha-\alpha_0 }\cdot \hat{C}M\eta ^{(k+1)}A_{k+1}\\
\end{aligned}
\end{equation*}
and
\begin{equation}\label{e.19}
\bar{a}\leq C_2\hat{C}MA_k.
\end{equation}
Take $\alpha_0=\alpha/2$ and then
\begin{equation}\label{e.20}
\begin{aligned}
\|v-\bar{a}x_n\|_{L^{\infty }(\Omega _{\eta^{k+1}})}=&\|v-\bar{a}x_n\|_{L^{\infty }(\Omega _{2\eta r})}\\
\leq & C_1\eta ^{\alpha_0 }\cdot \hat{C}M\eta ^{(k+1)}A_{k+1}+|\bar{a}r\omega(r)|\\
\leq &\left( C_1\eta ^{\alpha_{0} }+\frac{C_2\omega(\eta^{k})}{\eta^{1+\alpha_{0}}}\right)\cdot \hat{C}M\eta ^{(k+1)}A_{k+1}\\
\leq &\left( C_1\eta ^{\alpha_{0} }+\frac{C_2c_0}{\eta^{1+\alpha_{0}}}\right)\cdot \hat{C}M\eta ^{(k+1)}A_{k+1}.
\end{aligned}
\end{equation}

For $w$, by the Alexandrov-Bakel'man-Pucci maximum principle, we have
\begin{equation*}
  \begin{aligned}
\sup_{\Omega_{\eta^{k+1}}}w\leq \sup_{\tilde{\Omega} _{r}}w& \leq\|g\|_{L^{\infty }(\partial \Omega \cap \tilde{B}^{+}_{r})}+a_kr\omega(r)+C_3r\|f\|_{L^n(\tilde{\Omega}_{r})}\\
    &\leq Mr\omega_g(2r) +\sum_{i=0}^{k}|a_i-a_{i-1}|\eta^k\omega(\eta^k)+C_3Mr\omega_f(2r)\\
    &\leq M \eta^k \omega(\eta^k)+\bar{C}\hat{C}M\sum_{i=0}^{k}A_i\eta^k\omega( \eta^k)+C_3M\eta^k \omega(\eta^k).
  \end{aligned}
\end{equation*}
From the definition of $A_k$ (see \cref{e-Ak}), we have
\begin{equation*}
  \sum_{i=0}^{\infty} A_i\leq \sum_{i=1}^{\infty}\omega(\eta^i)+\eta^{\alpha_0}\sum_{i=0}^{\infty} A_i+c_0,
\end{equation*}
which implies
\begin{equation*}
\begin{aligned}
  \sum_{i=0}^{\infty} A_i\leq \frac{1}{1-\eta^{\alpha_0}}\sum_{i=1}^{\infty}\omega(\eta^i)+c_0
&=\frac{1}{1-\eta^{\alpha_0}}\sum_{i=1}^{\infty}\frac{\omega(\eta^i)
\left(\eta^{i-1}-\eta^i\right)}{\eta^{i-1}-\eta^i}+c_0\\
&=\frac{1}{\left(1-\eta^{\alpha_0}\right)\left(1-\eta\right)}\sum_{i=1}^{\infty}
\frac{\omega(\eta^i)\left(\eta^{i-1}-\eta^i\right)}{\eta^{i-1}}+c_0\\
&\leq \frac{1}{\left(1-\eta^{\alpha_0}\right)\left(1-\eta\right)}\int_{0}^{1}
\frac{\omega(r)dr}{r}+c_0\\
&\leq \frac{c_0}{\left(1-\eta^{\alpha_0}\right)\left(1-\eta\right)}+c_0\leq 3c_0,
\end{aligned}
\end{equation*}
provided
\begin{equation}\label{e-w}
\left(1-\eta^{\alpha_0}\right)\left(1-\eta\right)\geq 1/2.
\end{equation}
By the definition of $A_k$ again,
\begin{equation*}
  \omega(\eta^{k})\leq A_k\leq \frac{A_{k+1}}{\eta^{\alpha_0}}.
\end{equation*}
Hence,
\begin{equation}\label{e1.22}
  \begin{aligned}
\sup_{\Omega_{\eta^{k+1}}}w &\leq \frac{1}{\eta^{1+\alpha_0}}
M\eta^{k+1}A_{k+1}+\frac{3c_0\bar{C}}{\eta^{1+\alpha_0}}\hat{C}M\eta^{k+1}A_{k+1}+\frac{C_3 }{\eta^{1+\alpha_0}} M\eta^{k+1}A_{k+1}\\
&\leq \frac{C_3+1}{\eta^{1+\alpha_0}}
M\eta^{k+1}A_{k+1}+\frac{3c_0\bar{C}}{\eta^{1+\alpha_0}}\hat{C}M\eta^{k+1}A_{k+1}\\
&\leq\left(\frac{C_3+1}{\hat{C}\eta^{1+\alpha_0}}+\frac{3c_0\bar{C}}{\eta^{1+\alpha_0}}\right)
\hat{C}M\eta^{k+1}A_{k+1}.
\end{aligned}
\end{equation}

Let $\bar{C}=C_2$. Take $\eta $ small enough such that \cref{e-w} holds and
\begin{equation*}
  C_1\eta ^{\alpha_0 }\leq \frac{1}{4}.
\end{equation*}
Take $c_0$ small enough such that
\begin{equation*}
\frac{3c_0\bar{C}}{\eta^{1+\alpha_0}}\leq\frac{1}{4}.
\end{equation*}
Finally, take $\hat{C}$ large enough such that \cref{e.21} holds and
\begin{equation*}
  \frac{C_3+1}{\hat{C} \eta ^{1+\alpha_0}}\leq \frac{1}{4}.
\end{equation*}
Let $a_{k+1}=a_k+\bar{a}$. Then by combining \cref{e.20} and \cref{e1.22}, we have
\begin{equation*}
\begin{aligned}
u-a_{k+1}x_n&=u-a_kx_n-v+v-\bar{a}x_n=w+v-\bar{a}x_n\\
&\leq \hat{C}M\eta ^{(k+1)}A_{k+1} ~~\mbox{in}~~\Omega_{\eta^{k+1}}.
\end{aligned}
\end{equation*}
By induction, the proof is completed. \qed~\\

Next, we show that the exterior $C^{1,\mathrm{Dini}}$ condition is optimal, i.e., we give the~\\
\noindent\textbf{Proof of \Cref{t-3}.} For simplicity, we write $\omega=\omega_{\Omega}$ throughout this proof. Without loss of generality, we assume that $u\geq1$ in $B_{1/8}(e_n/2)$ and $\omega(1)\leq 1$. By \Cref{le-H-1} and noting that $u\geq 0$ on $T_1$, there exist universal constants $0<\delta_1<1$ and $c_1$ such that
\begin{equation}\label{e.3}
u(x)\geq c_1x_n ~\mbox{in}~B^{+}_{\delta_1}.
\end{equation}

To prove \cref{e-L-main}, we only need to prove the following: there exist universal constants $0<\eta \leq \delta_1$ and $0<c_0<1/2$, and a nonnegative sequence $\{a_k\}$ ($k\geq 0$) with $a_0=c_1$ such that for all $k\geq 1$,

\begin{equation}\label{L-e1.16}
u(x)\geq a_kx_n~\mbox{in}~B^+_{\eta^{k}}
\end{equation}
and
\begin{equation}\label{L-e1.17}
a_k=(1+c_0\omega(\eta^k))a_{k-1}.
\end{equation}

Indeed, for any $0<r<\eta$, there exists $k\geq 1$ such that $\eta^{k+1}\leq r \leq \eta^{k}$. Define
\begin{equation*}
\tilde{\omega}(r)=a_k=\prod_{i=1}^{k}(1+c_0\omega(\eta^i)).
\end{equation*}
Since $\omega$ don't satisfy the Dini condition, $\tilde{\omega}(r)\rightarrow +\infty$ as $r\rightarrow 0$. Then, for any $l\in\partial B_1$ with $l_n>0$,
\begin{equation*}
u(rl)\geq l_n\tilde{\omega}(r)r.
\end{equation*}
That is, \cref{e-L-main} holds.

Now, we prove \cref{L-e1.16} and \cref{L-e1.17} by induction. For $k=1$, by noting \cref{e.3}, they hold clearly. Suppose that they hold for $k$. We need to prove that they hold for $k+1$.

Let $r=\eta ^{k}$. By the weak Harnack inequality and \cref{L-e1.16},
\begin{equation*}\label{e.L.1}
u\geq c_2 a_kr~\mbox{on}~T_{r/2}\backslash T_{r/4}+\frac{r}{4}e_n,
\end{equation*}
where $c_2>0$ is a universal constant. Set
\begin{equation*}
Q=B_r\cap \left\{(x',x_n)\big| \frac{r}{4}<|x'|<\frac{r}{2},~-\frac{r}{4}\omega(\frac{r}{4})<x_n<\frac{r}{4}\right\}, T=\partial Q\cap \left\{x_n=\frac{r}{4}\right\}.
\end{equation*}
Then $Q\subset \Omega$. Let $v_1$ solve
\begin{equation*}
\left\{\begin{aligned}
 &M^{-}(D^2v_1,\lambda,\Lambda)=0 &&\mbox{in}~~Q; \\
 &v_1=c_2 a_kr &&\mbox{on}~~T;\\
 &v_1=0 &&\mbox{on}~~\partial Q\backslash T.
\end{aligned}
\right.
\end{equation*}
By the comparison principle, $u \geq v_1$ in $Q$ clearly. From \Cref{le-H-2-2},
\begin{equation*}
v_1(x)\geq c_3a_k(x_n+\frac{r}{4}\omega(\frac{r}{4})) ~\mbox{in}~\tilde{Q},
\end{equation*}
where $c_3$ is a universal constant and
\begin{equation*}
\tilde{Q}=Q\cap \left\{(x',x_n)\big| \frac{r}{3}<|x'|<\frac{5r}{12},~-\frac{r}{4}\omega(\frac{r}{4})<x_n\leq -\frac{r}{4}\omega(\frac{r}{4})+\delta_2r\right\}.
\end{equation*}
Here, $\delta_2<1/8$ is a universal constant. Without loss of generality, we assume that $\omega(1)\leq \delta_2$. Hence,
\begin{equation*}
  u\geq v_1\geq \frac{c_3a_kr}{4}\omega(\frac{r}{4})~\mbox{on}~T_{5r/12}\backslash T_{r/3}.
\end{equation*}

Next, let $v_2$ solve
\begin{equation*}
\left\{\begin{aligned}
 &M^{-}(D^2v_2,\lambda,\Lambda)=0 &&\mbox{in}~~B_r^+; \\
 &v_2=c_3a_k\frac{r}{4}\omega(\frac{r}{4}) &&\mbox{on}~~T_{5r/12}\backslash T_{r/3};\\
 &v_2=0 &&\mbox{on}~~\partial B_r^+\backslash (T_{5r/12}\backslash T_{r/3}).
\end{aligned}
\right.
\end{equation*}
By the comparison principle again,
\begin{equation*}
u\geq a_kx_n+v_2 ~\mbox{in}~ B_r^+.
\end{equation*}
From \Cref{le-H-2},
\begin{equation*}
v_2(x)\geq c_0a_k\omega(\frac{r}{4})x_n~\mbox{in}~B^+_{\delta_2r}.
\end{equation*}

Take $\eta=\min (\delta_1,\delta_3,1/4)$ and then
\begin{equation*}
u(x)\geq  a_kx_n+v_2(x)\geq a_kx_n+c_0a_k\omega(\frac{r}{4})x_n
\geq \left(1+c_0\omega(\eta^{k+1})\right)a_kx_n~\mbox{in}~B^+_{\eta^{k+1}}.
\end{equation*}
Thus, \cref{L-e1.16} and \cref{L-e1.17} hold for $k+1$. By induction, the proof is completed. \qed~\\

\section{Hopf lemma}
The proof of the Hopf lemma is similar to that of the boundary Lipschitz regularity. Here, we focus on the curved boundary toward the interior of the domain. Now, we give the~\\
\noindent\textbf{Proof of \Cref{t-H-2}.} For simplicity, we write $\omega=\omega_{\Omega}$ throughout this proof. As before, from the Dini condition \cref{e-dini}, there exists $r_1>0$ such that
\begin{equation}\label{e-H-dini-2}
 \omega(r_1)\leq c_0 ~~\mbox{and}~~ \int_{0}^{r_1}\frac{\omega(r)}{r}dr\leq c_0,
\end{equation}
where $c_0\leq 1/4$ is a universal small constant to be specified later. By a proper scaling, we assume that $r_1=1$ and $u\geq1$ in $B_{1/8}(e_n/2)$ without loss of generality.

Let $\Omega^+ _{r}=\Omega \cap B^+_{r}$. To prove \cref{e-H-main}, we only need to prove the following: there exist universal constants $0<\alpha_{0}, \eta < 1$, $\bar{C}$ and $\tilde{a}>0$,  and a nonnegative sequence $\{a_k\}$ ($k\geq -1$) such that for all $k\geq 0$,

\begin{equation}\label{e-H1.16}
\inf_{\Omega ^+_{\eta^{k+1}}}(u-\tilde{a}x_n+a_kx_n)\geq -\eta ^{k}A_k,
\end{equation}
\begin{equation}\label{e-H1.17}
a_k-a_{k-1}\leq \bar{C}A_k
\end{equation}
and
\begin{equation}\label{e-H-atilde}
a_k\leq \frac{\tilde{a}}{2},
\end{equation}
where $A_k$ is defined as in \cref{e-Ak}.

Indeed, for any $l\in\partial B_1$ with $l_n>0$, there exists $k_0\geq 1$ such that for any $k\geq k_0$,
\begin{equation*}
\frac{A_k}{l_n\eta}\leq \frac{\tilde{a}}{4}.
\end{equation*}
Take $\delta=\eta^{k_0}$. For $0<r<\delta$, there exists $k\geq k_0$ such that $  \eta^{k+1}\leq r \leq \eta^{k}$. Then by \cref{e-H1.16},
\begin{equation*}
\begin{aligned}
u(rl)\geq \tilde{a}l_nr-a_kl_nr-\eta^kA_k
\geq \frac{\tilde{a}l_nr}{2}-\frac{A_kr}{\eta}
\geq \frac{\tilde{a}l_n}{4} r.
\end{aligned}
\end{equation*}
That is, \cref{e-H-main} holds.

Now, we prove \crefrange{e-H1.16}{e-H-atilde} by induction. In the following proof, $C,C_1,C_2...,$ always denote universal constants.
%
From \Cref{le-H-1}, there exist universal constants $\delta_1>0$ and $0<c_1<1/2$ such that
\begin{equation*}
u(x)\geq c_1(x_n-\omega(1)) ~~\mbox{in}~~B_{\delta_1}^++\omega(1)e_n.
\end{equation*}
Note that $u\geq 0$ and hence
\begin{equation*}
  u(x)\geq c_1(x_n-\omega(1)) ~~\mbox{in}~~\Omega_{\delta_1}^+.
\end{equation*}

Set $\tilde{a}=c_1$, $a_{-1}=a_0=0$ and
\begin{equation}\label{e-H-1}
\eta\leq \delta_1.
\end{equation}
Then \crefrange{e-H1.16}{e-H-atilde} hold for $k=0$. Suppose that they hold for $k$. We need to prove that they hold for $k+1$.

Let $r=\eta ^{k+1}$ and $v$ solve
\begin{equation*}
\left\{\begin{aligned}
 &M^{-}(D^2v,\lambda,\Lambda)=0 &&\mbox{in}~~B^{+}_{r}; \\
 &v=0 &&\mbox{on}~~T_{r};\\
 &v=-\eta ^{k}A_k &&\mbox{on}~~\partial B^{+}_{r}\backslash T_{r}.
\end{aligned}
\right.
\end{equation*}
Let $w=u-\tilde{a}x_n+a_kx_n-v$. Then $w$ satisfies (note that $u\geq 0$ and $v\leq 0$)
\begin{equation*}
    \left\{
    \begin{aligned}
      &M^{-}(D^2w,\lambda/n,\Lambda)\leq 0 &&\mbox{in}~~ \Omega^{+}_{r}; \\
      &w\geq -\tilde{a}x_n+a_kx_n &&\mbox{on}~~\partial \Omega \cap B^{+}_{r};\\
      &w\geq 0 &&\mbox{on}~~\partial B^{+}_{r}\cap \bar{\Omega}.
    \end{aligned}
    \right.
\end{equation*}

In the following arguments, we estimate $v$ and $w$ respectively. By the boundary $C^{1,\alpha}$ estimate for $v$ (see \Cref{l-1}) and the maximum principle, there exist $0<\alpha<1$ (depending only on $n,\lambda$ and $\Lambda$) and $\bar{a}\geq 0$  such that (note that $A_k\leq A_{k+1}/\eta^{\alpha_0}$)
\begin{equation*}
\begin{aligned}
    \|v+\bar{a}x_n\|_{L^{\infty }(B^+ _{\eta r})}\leq C \eta ^{1+ \alpha}\|v\|_{L^{\infty }(B^{+}_{r})}
    \leq C\eta ^{\alpha}\cdot \eta ^{(k+1)}A_k
    \leq C_1\eta ^{\alpha-\alpha_0 }\cdot \eta ^{(k+1)}A_{k+1}\\
\end{aligned}
\end{equation*}
and
\begin{equation}\label{e.H-19}
\bar{a}\leq C_2A_k/\eta.
\end{equation}
Take $\alpha_0=\alpha/2$ and then
\begin{equation}\label{e.H-20}
\begin{aligned}
\|v+\bar{a}x_n\|_{L^{\infty }(\Omega^+_{\eta^{k+2}})}\leq\|v+\bar{a}x_n\|_{L^{\infty }(B^+_{\eta r})}\leq C_1\eta ^{\alpha_0 }\cdot\eta ^{(k+1)}A_{k+1}.
\end{aligned}
\end{equation}

For $w$, by the maximum principle, we have
\begin{equation*}
  \begin{aligned}
\inf_{\Omega^+_{\eta^{k+2}}}w&\geq \inf_{\Omega^+ _{r}}w \geq -\tilde{a}r\omega(r)= -\tilde{a}\eta^{k+1}\omega(\eta^{k+1}).
  \end{aligned}
\end{equation*}
As before, from the definition of $A_k$ (see \cref{e-Ak}),
\begin{equation*}
  \omega(\eta^{k+1})\leq A_{k+1}.
\end{equation*}
Hence,
\begin{equation}\label{e-H-1.22}
  \begin{aligned}
\inf_{\Omega^+_{\eta^{k+2}}}w \geq -\tilde{a}\eta^{k+1}A_{k+1}
\geq -\frac{1}{2}\eta^{k+1}A_{k+1}.
\end{aligned}
\end{equation}

Take $\eta $ small enough such that \cref{e-H-1} holds,
\begin{equation*}
  C_1\eta ^{\alpha_0 }\leq 1/2
\end{equation*}
and
\begin{equation}\label{e-H-2}
\left(1-\eta^{\alpha_0}\right)\left(1-\eta\right)\geq 1/2.
\end{equation}
Let $\bar{C}=C_2/\eta$. As before, by noting that \cref{e-H-2}, we have
\begin{equation*}
  \sum_{i=0}^{\infty} A_i\leq 3c_0.
\end{equation*}
Thus,
\begin{equation*}
  a_k\leq \sum_{i=0}^{k}|a_i-a_{i-1}|\leq \bar{C}\sum_{i=0}^{\infty} A_i\leq 3c_0 \bar{C}.
\end{equation*}
Finally, take $c_0$ small enough such that
\begin{equation*}
3c_0 \bar{C}\leq\frac{\tilde{a}}{2}.
\end{equation*}

Let $a_{k+1}=a_k+\bar{a}$. Then by combining \cref{e.H-20} and \cref{e-H-1.22}, we have
\begin{equation*}
\begin{aligned}
u-\tilde{a}x_n+a_{k+1}x_n&=u-\tilde{a}x_n+a_kx_n-v+v+\bar{a}x_n=w+(v+\bar{a}x_n)\\
&\geq -\eta ^{(k+1)}A_{k+1} ~~\mbox{in}~~\Omega^+_{\eta^{k+2}}.
\end{aligned}
\end{equation*}
By induction, the proof is completed. \qed~\\

Next, we show that the interior $C^{1,\mathrm{Dini}}$ condition is optimal and give the~\\
\noindent\textbf{Proof of \Cref{t-H-3}.} For simplicity, we write $\omega=\omega_{\Omega}$ throughout this proof. Without loss of generality, we assume that $\omega(1)\leq 1$ and $u\leq1$ in $\Omega\cap B_1$. Let $v_1$ solve
\begin{equation*}
\left\{\begin{aligned}
 &M^{+}(D^2v_1,\lambda,\Lambda)=0 &&\mbox{in}~~B^{+}_{1}; \\
 &v_1=0 &&\mbox{on}~~T_{1};\\
 &v_1=1 &&\mbox{on}~~\partial B^{+}_{1}\backslash T_{1}.
\end{aligned}
\right.
\end{equation*}
By \Cref{l-1}, there exist universal constants $0<\delta_1<1$ and $C_1$ such that
\begin{equation*}
v_1(x)\leq C_1x_n ~\mbox{in}~B^{+}_{\delta_1}.
\end{equation*}
Hence,
\begin{equation*}
  u(x)\leq v_1(x)\leq C_1x_n ~\mbox{in}~\Omega\cap B^{+}_{\delta_1}.
\end{equation*}

Let $\Omega^+ _{r}=\Omega \cap B^+_{r}$. To prove \cref{e-H-main-2}, we only need to prove the following: there exist universal constants $0<\eta \leq \delta_1$ and $0<c_0<1/2$, and a nonnegative sequence $\{a_k\}$ ($k\geq 0$) such that for all $k\geq 1$,

\begin{equation}\label{H-e1.16}
u(x)\leq a_kx_n~\mbox{in}~\Omega^+_{\eta^{k}}
\end{equation}
and
\begin{equation}\label{H-e1.17}
a_k= (1-c_0\omega(\eta^k))a_{k-1}.
\end{equation}

Indeed, for any $0<r<\eta$, there exists $k\geq 1$ such that $\eta^{k+1}\leq r \leq \eta^{k}$. Define
\begin{equation*}
\tilde{\omega}(r)=a_k=a_0\prod_{i=1}^{k}(1-c_0\omega(\eta^i)).
\end{equation*}
Since $\omega$ don't satisfy the Dini condition, $\tilde{\omega}(r)\rightarrow 0$ as $r\rightarrow 0$. Hence,
\begin{equation*}
\sup_{\Omega_{r}} u\leq r\tilde\omega(r), ~~\forall~ 0<r<\eta.
\end{equation*}
That is, \cref{e-H-main-2} holds.

Now, we prove \cref{H-e1.16} and \cref{H-e1.17} by induction. For $k=0$, by setting $a_0=C_1$ and $a_{-1}=2C_1$, they hold clearly. Suppose that they hold for $k$. We need to prove that they hold for $k+1$.

Let $r=\eta ^{k}$,
\begin{equation*}
B^+=B_r^+\cap \left\{(x',x_n)\big| |x'|<\frac{r}{4}~\mbox{or}~x_n>\frac{r}{4}\omega(\frac{r}{4})\right\}, T=T_{r}\backslash T_{r/4}+\frac{r}{4}\omega(\frac{r}{4})e_n
\end{equation*}
and $v$ solve
\begin{equation*}
\left\{\begin{aligned}
 &M^{-}(D^2v,\lambda,\Lambda)=0 &&\mbox{in}~~B^+; \\
 &v=\frac{a_kr}{4}\omega(\frac{r}{4}) &&\mbox{on}~~T;\\
 &v=0 &&\mbox{on}~~\partial B^+\backslash T.
\end{aligned}
\right.
\end{equation*}
Take $w=a_kx_n-u-v$ and $w$ satisfies
\begin{equation*}
    \left\{
    \begin{aligned}
      &w\in \bar{S}(\lambda,\Lambda, 0) &&\mbox{in}~~ \Omega \cap B^{+}; \\
      &w\geq 0 &&\mbox{on}~~\partial \left(\Omega \cap B^{+}\right).
    \end{aligned}
    \right.
\end{equation*}
Indeed, by \cref{H-e1.16} and $v= 0$ on $\partial B^{+}\cap \Omega \cap \left\{x_n>r\omega(r/4)/4\right\}$, we have $w\geq 0$ on $\partial B^{+}\cap \Omega \cap \left\{x_n>r\omega(r/4)/4\right\}$. In addition, since $u=0$ on $\partial \Omega\cap B^+$ and $v\leq a_kx_n$ in $B^+$, we have $w\geq 0$ on $\partial \Omega \cap B^{+}$. Hence, by the maximum principle, $w\geq 0$ in $\Omega \cap B^{+}$.

For $v$, by \Cref{le-H-2}, there exist universal constants $\delta_2$ and $c_0$ such that
\begin{equation*}
v\geq c_0a_k\omega(\frac{r}{4})x_n ~\mbox{in}~B^+_{\delta_2r}.
\end{equation*}
Take $\eta=\min (\delta_1,\delta_2,1/4)$ and then
\begin{equation*}
u(x)\leq  a_kx_n-v(x)\leq a_kx_n-c_0a_k\omega(\frac{r}{4})x_n
\leq \left(1-c_0\omega(\eta^{k+1})\right)a_kx_n~\mbox{in}~\Omega^+_{\eta^{k+1}}.
\end{equation*}
Thus, \cref{H-e1.16} and \cref{H-e1.17} hold for $k+1$. By induction, the proof is completed. \qed~\\

\section*{Declarations}
\noindent \textbf{Ethical Approval } not applicable.
~\\

\noindent \textbf{Competing interests} The authors declare that the authors have no competing interests as defined by Springer, or other interests that might be perceived to influence the results and/or discussion reported in this paper.
~\\

\noindent \textbf{Authors' contributions} Lian and Zhang discussed the problem many times and the main idea originated from their discussions. Lian wrote the main manuscript text and Zhang made some revisions.
~\\

\noindent \textbf{Funding} This research is supported by the China Postdoctoral Science Foundation (Grant No. 2021M692086 and 2022M712081), the National Natural Science Foundation of China (Grant No. 12031012, 11831003 and 12171299) and the Institute of Modern Analysis-A Frontier Research Center of Shanghai.
~\\

\noindent \textbf{Availability of data and materials} not applicable.

\bibliographystyle{amsplain}
\bibliography{L-H}

\end{document}